\documentclass[11pt, a4paper]{article}
\usepackage{amsthm}
\usepackage{mathrsfs}
\usepackage{amsmath,amssymb,latexsym,color}
\usepackage[colorlinks,
linkcolor=blue,
anchorcolor=green,
citecolor=magenta
]{hyperref}
\usepackage{graphicx}
\usepackage{tikz}
\usetikzlibrary{calc}
\usepackage{CJK}

\usepackage{authblk}

\long\def\delete#1{}

\usepackage{color}

\definecolor{Blue}{rgb}{0,0,1}
\definecolor{Red}{rgb}{1,0,0}
\definecolor{DarkGreen}{rgb}{0,0.6,0}
\definecolor{DarkYellow}{rgb}{1,1,0.2}
\definecolor{DarkPurple}{rgb}{.6,0,1}

\usepackage{xcolor}
\usepackage[normalem]{ulem}

\usepackage{enumerate}
\usepackage{cleveref}
\crefformat{section}{\S#2#1#3}
\crefformat{subsection}{\S#2#1#3}
\crefformat{subsubsection}{\S#2#1#3}
\crefrangeformat{section}{\S\S#3#1#4 to~#5#2#6}
\crefmultiformat{section}{\S\S#2#1#3}{ and~#2#1#3}{, #2#1#3}{ and~#2#1#3}
\def\ma{\mathcal{A}}
\def\mb{\mathcal{B}}

\def\md{\mathcal{D}}
\def\me{\mathcal{E}}
\def\mf{\mathcal{F}}
\def\mg{\mathcal{G}}
\def\mh{\mathcal{H}}
\def\mi{\mathcal{I}}

\def\mm{\mathcal{M}}

\def\ms{\mathcal{S}}
\def\mt{\mathcal{T}}
\def\mv{\mathcal{V}}
\def\mw{\mathcal{W}}

\def\my{\mathcal{Y}}

\def\bs{\setminus}
\def\mmu{\mathcal{U}}

\numberwithin{equation}{section}

\newtheorem{thm}{Theorem}[section]
\newtheorem{cl}[thm]{Claim}
\newtheorem{lem}[thm]{Lemma}

\newtheorem{cor}[thm]{Corollary}
\newtheorem{pr1}[thm]{Proposition}
\newtheorem{ex}[thm]{\bf Example}
\theoremstyle{remark}
\newtheorem{re}{\bf Remark}

\begin{document}
	\setcounter{page}{1}
	\renewcommand{\thefootnote}{}
	\newcommand{\remark}{\vspace{2ex}\noindent{\bf Remark.\quad}}
	\renewcommand{\abovewithdelims}[2]{%
		\genfrac{[}{]}{0pt}{}{#1}{#2}}

	%-------------------  First Head  -----------------------------------------
	
	\def\qed{\hfill$\Box$\vspace{11pt}}
	
	\title {\bf $s$-almost $t$-intersecting families for finite sets}

	\author{Dehai Liu\thanks{E-mail: \texttt{liudehai@mail.bnu.edu.cn}}\   \textsuperscript{a}}
	\author{Kaishun Wang\thanks{ E-mail: \texttt{wangks@bnu.edu.cn}}\ \textsuperscript{a}}
	\author{Tian Yao\thanks{Corresponding author. E-mail: \texttt{tyao@hist.edu.cn}}\ \textsuperscript{b}}
	\affil{ \textsuperscript{a} Laboratory of Mathematics and Complex Systems (Ministry of Education), School of
		Mathematical Sciences, Beijing Normal University, Beijing 100875, China}
	
	\affil{ \textsuperscript{b} School of Mathematical Sciences, Henan Institute of Science and Technology, Xinxiang 453003, China}
	\date{}
	
	\openup 0.5\jot
	\maketitle

\begin{abstract}

	A family $\mathcal{F}$ of $k$-subsets of an $n$-set is called $s$-almost $t$-intersecting if each member is $t$-disjoint with at most $s$ members. 	
	In this paper, we prove that, if $\left|\mathcal{F}\right|$ is maximum, then $\mathcal{F}$ consists of all $k$-subsets containing a fixed $t$-subset.  	
	Consequently, it is natural to consider the maximum-sized $\mathcal{F}$ with $\left|\bigcap_{F\in\mathcal{F}} F\right|<t$. The famous Hilton-Milner theorem settles the case where $\mathcal{F}$ is $t$-intersecting. We characterize the remaining case completely.

	\vspace{2mm}
	
	\
	
	
\end{abstract}

\section{Introduction}

Let $n$, $k$, $t$ and $s$ be positive integers with $n\geq  k\geq t+1$.  Write $[n]=\left\{1,2,\ldots,n\right\}$.  For a set $X$, the family of all $k$-subsets of $X$ is denoted by $\binom{X}{k}$.  Let $K_n$ denote the complete graph with $n$ vertices, $S_n$ the star graph with $n$ vertices, and $K_{s,t}$ the complete bipartite graph with parts of sizes $s$ and $t$.

The \textit{generalized Kneser graph} $K(n,k,t)$ is a graph whose vertex set is $\binom{[n]}{k}$, where two vertices $A$ and $B$ are adjacent if and only if $\left|A\cap B\right|<t$. 
For a family  $\mf\subseteq \binom{[n]}{k}$, write $K(n,k,t)[\mf]$ for  the subgraph of $K(n,k,t)$ induced  by  $\mf$. 
Extensive results   have been achieved concerning the largest family $\mf$  under the condition that $K(n,k,t)[\mf]$ has no  subgraph isomorphic to a  graph $G$.  
The case  $G=K_{2}=S_{2}$ corresponds to the famous Erd\H{o}s-Ko-Rado theorem \cite{MR0140419, MR0519277,MR0771733} and the  Ahlswede-Khachatrian Complete Intersection theorem \cite{2406295}.
The case $t=1$ and $G=K_{r}$ is the well-known Erd\H{o}s matching conjecture \cite{2410242}, which has been proved in  \cite{24121102,2501163,2412114,2501161} for sufficiently large $n$. 
The case $t=1$ and $G=K_{p,q}$ has also been studied in \cite{2410232,2412111}  for  Erd\H{o}s-Ko-Rado type results.

In this paper, we focus on the case that $G$ is a star graph.  
We say a family $\mf\subseteq \binom{[n]}{k}$ is \textit{$s$-almost $t$-intersecting} if $\left|\left\{H\in\mf: \left|H\cap F\right|<t\right\}\right|\leq s$ for any $F\in\mf$, which holds   if and only if $K(n,k,t)[\mf]$ is $S_{s+2}$-free. A family $\mf\subseteq \binom{[n]}{k}$ is said to be \textit{$t$-intersecting} if $K(n,k,t)[\mf]$ is $S_{2}$-free. Consequently, each $t$-intersecting family is $s$-almost $t$-intersecting for all $s$.

In \cite{2406095}, Gerbner et al.  characterized the maximum-sized $s$-almost $t$-intersecting family  for $t=1$. 
Our first result extends their result to $t\geq 1$. For convenience, set
\begin{equation*}
	\begin{aligned}
		L_{1}(k,t,s)&=(k-t+1)((t+1)+3s^{\frac{1}{k-t}}),\\
		L_{2}(k,t,s)&=2(k-t+1)^{2}((t+1)^{2}+3s^{\frac{1}{k-t-1}}/ (k-t+1)).
	\end{aligned}
\end{equation*}
\begin{thm}\label{2509211}
	Let $n$, $k$, $t$ and $s$ be positive integers with $k\geq t+1$ and $n\geq L_{1}(k,t,s)$. If $\mf\subseteq \binom{[n]}{k}$ is an $s$-almost $t$-intersecting family with the maximum size, then there exists $X\in\binom{[n]}{t}$ such that $\mf=\left\{F\in\binom{[n]}{k}: X\subseteq F\right\}$.
\end{thm}

Motivated by  Theorem \ref{2509211}, it is natural to consider   the $s$-almost $t$-intersecting family $\mf\subseteq \binom{[n]}{k}$ with $\left|\bigcap_{F\in\mf}F\right|<t$.   Note that, under the condition that $\mf$ is $t$-intersecting, the extremal structure has been characterized by the Hilton-Milner theorem \cite{2410236, 2407161,2407162}. We therefore turn our attention to the remaining situation.

Suppose $t=1$. For $k=2$ and $s\geq 14$, the extremal family was  determined by Scott and Wilmer \cite{2412121}. Subsequently, Alishahi and Taherkhani \cite{2410232}  provided an asymptotically optimal upper bound on $\left|\mf \right|$. For $k\geq 3$ and $s=1$, Frankl and  Kupavskii \cite{2406096} described  the maximum-sized $\mf$. A significant step was taken in  \cite{2410233} by Taherkhani, who gave  the extremal structure of $\mf$  for  $k\geq3$ and sufficiently large $n$.

  Our following results    characterize the extremal structure of $s$-almost $t$-intersecting families that are not $t$-intersecting. To present them, we introduce some notations.  Suppose that $X$, $Y$ and $Z$ are subsets of $[n]$. Write
\begin{equation*}
	\begin{aligned}
		\mh(X,Y;k,t)&=\left\{F\in\binom{[n]}{k}: Y\subseteq F,\ \left|F\cap X\right|\geq t+1\right\},\\
		\mm_{1}(X,Y,Z;t)&=\left\{ F\in\binom{X}{t+1}: Z\subseteq F,\ \left|F\cap Y\right|\geq t\right\},\\
		\mm_{2}(X,Y,Z; t)&=\left\{F\in \binom{X}{t+1}: Z\subseteq F,\ Y\nsubseteq F \right\}. 
	\end{aligned}
\end{equation*}

\begin{thm}\label{2509212}
	Let $n$, $k$, $t$ and $s$ be positive integers with $k\geq t+2$ and $n\geq L_{2}(k,t,s)$. Suppose that $\mf\subseteq \binom{[n]}{k}$ is an $s$-almost $t$-intersecting family but not $t$-intersecting. If $\mf$ has the maximum size, then there exist $X\in\binom{[n]}{k+1}$ and $Y\in\binom{X}{t}$ such that	$$\mf=\mh(X,Y;k,t)\cup\ma\cup\mb,$$ 
	where $\ma\subseteq\left\{ F\in\binom{[n]}{k}: F\cap X=Y\right\}$ with $\left|\ma\right|=s$ and $\mb\subseteq \left\{ F\in\binom{X}{k}: Y \nsubseteq F\right\}$ with  $\left|\mb\right|=\min\left\{t,s\right\}$.
\end{thm}

\begin{thm}\label{2509213}
Let $n$, $t$ and $s$ be positive integers with $s\geq 7$ and $n\geq t+3$. If $\mf\subseteq \binom{[n]}{t+1}$ is an $s$-almost $t$-intersecting family but not $t$-intersecting, then $\left|\mf\right|\leq 2s+3$. Moreover, equality holds if and only if one of the following holds.
\begin{enumerate}[\normalfont(i)]
	\item $n\geq t+s+2$ and  $\mf=\mm_{1}(X,Y,Z;t)$ for some $X\in \binom{[n]}{t+s+2}$, $Y\in\binom{X}{t+1}$ and $Z\in\binom{Y}{t-1}$.
	\item $t\geq s$ and  $\mf=\mm_{2}(X,Y,Z;t)$ for some $X\in \binom{[n]}{t+3}$, $Y\in\binom{X}{2}$ and $Z\in\binom{X\bs Y}{t-s}$.
\end{enumerate}
\end{thm}

The rest of this paper is organized as follows.  In Section \ref{2510201}, we show some  results which are required in the proofs.   Theorems \ref{2509211}, \ref{2509212} and \ref{2509213}  are proved in Sections \ref{2510203}, \ref{2510204} and  \ref{2510205} respectively.    To ensure a smoother reading experience, some   inequalities needed in this paper are left in Section \ref{2510202}.

\section{Preliminaries}\label{2510201}

In this section, we review some earlier results concerning extremal set theory, and prove some auxiliary results which are used in  the present paper. 
For a family $\mf\subseteq \binom{[n]}{k}$  and a set $G\in\binom{[n]}{k}$, write
$$\mi_{\mf}(t)=\left\{F\in\mf: \left|F\cap H\right|\geq t\ \textnormal{for any $H\in\mf$}\right\},$$
$$\md_{\mf}(t)=\left\{F\in\mf: \left|F\cap H\right|< t\ \textnormal{for some $H\in\mf$}\right\},$$
$$\md_{\mf}(G;t)=\left\{F\in\mf: \left|F\cap G\right|<t\right\}.$$
Note that  $\mf$ is the disjoint union of $\mi_{\mf}(t)$ and $\md_{\mf}(t)$, and  $\md_{\mf}(F;t)\subseteq \md_{\mf}(t)$ for any $F\in\mf$.

One  result we need  concerns the $t$-intersecting family $\mf \subseteq \binom{[n]}{k}$   with $\left|\bigcap_{F\in\mf}F\right|<t$. For the sake of simplicity, we set
\begin{align}
g_{1}(n,k,t)&=	(t+2)\binom{n-t-1}{k-t-1}-(t+1)\binom{n-t-2}{k-t-2},\label{2510162}\\ 
g_{2}(n,k,t)&=	\binom{n-t}{k-t}-\binom{n-k-1}{k-t}+t.\label{2510163}
\end{align}

\begin{thm}\label{2510151} \textnormal{(\cite{2410236, 2407161,2407162})}
Let $n$, $k$ and $t$ be positive integers with $k\geq t+1$ and $n\geq (t+1)(k-t+1)$. If $\mf \subseteq \binom{[n]}{k}$ is a $t$-intersecting family with $\left|\bigcap_{F\in\mf}F\right|<t$, then 
 $$\left|\mf\right|\leq \max\left\{g_{1}(n,k,t), g_{2}(n,k,t)\right\}.$$
\end{thm}

Another ingredient we need is the skew version of Bollob\'{a}s theorem \cite{2510153}.  Applying this theorem, we  deduce an upper bound on $\left|\md_{\mf}(t)\right|$  of an $s$-almost $t$-intersecting family $\mf$.

\begin{thm}\label{2510152} \textnormal{(\cite{2510153})}
	Let $n$, $k$, $t$ and $m$ be positive integers with $n\geq k\geq t$. Suppose that $A_{1}, A_{2},\ldots, A_{m}$ and $B_{1}, B_{2},\ldots, B_{m}$ are $k$-subsets of $[n]$. If  $\left|A_{i}\cap B_{i}\right|< t$ for any $1\leq i\leq m$, and $\left|A_{i}\cap B_{j}\right|\geq t$ for any $1\leq j<i\leq m$, then $m\leq \binom{2k-2t+2}{k-t+1}$. 
\end{thm}

\begin{lem}\label{2509214}
	Let $n$, $k$, $t$ and $s$ be positive integers with $n\geq k\geq t$. If  $\mf\subseteq \binom{[n]}{k}$ is an $s$-almost $t$-intersecting family, then $\left|\md_{\mf}(t)\right|\leq s\binom{2k-2t+2}{k-t+1}$. 
\end{lem}
\begin{proof}
	If $\md_{\mf}(t)=\emptyset$, then the required result is clear. Next assume $\md_{\mf}(t)\neq \emptyset$.
	Set $\me_{1}=\md_{\mf}(t)$.  Choose $A_{i}$, $B_{i}$ and $\me_{i+1}$ by repeating the following steps:
	$$A_{i}\in \me_{i},\ B_{i}\in\md_{\mf}(A_{i};t),\ \me_{i+1}=\me_{i}\bs\md_{\mf}(B_{i};t).$$
	Since $\left| \me_{i+1}\right|< \left| \me_{i}\right|$, there exists a positive integer $m$ such that $\me_{m+1}=\emptyset$.	Finally, we  get two sequences of subsets $A_{1},A_{2},\ldots,A_{m}$ and $B_{1}, B_{2},\ldots,B_{m}$ such that
	\begin{enumerate}[\normalfont(i)]
		\item $\left|A_{i}\cap B_{i} \right|<t$ for any $i\in[m]$.
		\item $\left|A_{i}\cap B_{j}\right|\geq t$ for any $i,j\in[m]$ with $j<i$.
		\item $\md_{\mf}(t)=\cup _{i=1}^{m}\md_{\mf}(B_{i};t)$.
	\end{enumerate}
	
	According to Theorem \ref{2510152}, two sequences of subsets $A_{1},A_{2},\ldots,A_{m}$ and $B_{1}, B_{2},\ldots,B_{m}$ satisfying  (i) and (ii) must have $m\leq \binom{2k-2t+2}{k-t+1}$.	Notice that $\mf$ is $s$-almost $t$-intersecting. 
	From (iii), we obtain 
	$\left|\md_{\mf}(t)\right|\leq \sum_{i=1}^{m}\left|\md_{\mf}(B_{i};t)\right|\leq s \binom{2k-2t+2}{k-t+1}$,
	as desired.
\end{proof}

For small  $k$, we proceed by showing more precise upper bounds. Suppose that $A$ and $B$ are two subsets of $[n]$. Write
\begin{equation}\label{2512291}
	\mw(A,B)=\left\{ \left(A\cap B\right)\cup\left\{ a,b\right\}: a\in A\bs B,\ b\in B\bs A\right\}.
\end{equation}
Observe that if  $\left|A\cap B\right|=t-1$, then 
\begin{equation}\label{2510021}
	\mw(A,B)=\left\{F\in\binom{[n]}{t+1}: \left|F\cap A\right|\geq t,\ \left|F\cap B\right|\geq t\right\}.
\end{equation}

\begin{lem}\label{2509215}
	Let $n$, $k$, $t$ and $s$ be positive integers with $ k\geq t$ and $n\geq 2k-t+1$. Suppose that $\mf\subseteq \binom{[n]}{k}$ is an $s$-almost $t$-intersecting family but not $t$-intersecting. The following hold.
	\begin{enumerate}[\normalfont(i)]
		\item If $k=t+1$, then $\left|\mf\right|\leq 2s+4$.
		\item If $k=t+2$, then $\left|\md_{\mf}(t)\right|\leq 11s+9t+18$.
	\end{enumerate}
\end{lem}
\begin{proof}
	 Pick $A,B\in\md_{\mf}(t)$  with $\left|A\cap B\right|=\min\left\{ \left|F_{1}\cap F_{2}\right|: F_{1}, F_{2}\in\md_{\mf}(t)\right\}$, which implies $\left|A\cap B\right|\leq t-1$. For each $T\subseteq [n]$ with $\left|T\cap A\right|\geq t$ and $\left|T\cap B\right|\geq t$, we have 
	\begin{equation}\label{2509247}
		\left|T\right|\geq \left|T\cap \left(A\cup B\right)\right|=\left|T\cap A\right|+\left|T\cap B\right|-\left|T\cap A\cap B\right|\geq 2t-\left| A\cap B\right|.
	\end{equation}

 (i) If $\left|A\cap B\right|\leq t-2$, then by (\ref{2509247}), we know $\mf\subseteq \md_{\mf}(A;t)\cup \md_{\mf}(B;t)$. Hence $\left|\mf\right|\leq 2s$.	If $\left|A\cap B\right|=t-1$, by (\ref{2510021}), we get
	$\mf\subseteq \mw(A,B) \cup \md_{\mf}(A;t)\cup \md_{\mf}(B;t)$. Thus
	 $\left|\mf\right|\leq 2s+4$. 
	
	(ii) If $ \md_{\mf}(t)\subseteq \md_{\mf}(A;t)\cup \md_{\mf}(B;t)$, then $\left|\md_{\mf}(t)\right|\leq 2s$, as desired. Next, we consider $\md_{\mf}(t)\bs \left(\md_{\mf}(A;t)\cup \md_{\mf}(B;t)\right)\neq \emptyset$. This together with  (\ref{2509247}) yields $\left|A\cap B\right|\geq t-2$. 
	
	\medskip
	\noindent \textbf{Case 1.} $\left|A\cap B\right|=t-2$.
	\medskip
	
	Pick    $D\in\md_{\mf}(t)\bs\left(\md_{\mf}(A;t)\cup\md_{\mf}(B;t)\right)$. By (\ref{2509247}), we have  $\left|D\cap A\right|=\left|D\cap B\right|=t$ and $\left|D\cap A\cap B\right|=t-2$. The latter implies $A\cap B\subseteq D$. Hence
	$$\md_{\mf}(t)\bs\left(\md_{\mf}(A;t)\cup\md_{\mf}(B;t)\right)\subseteq \left\{ \left(A\cap B\right)\cup D_{1}\cup D_{2}: D_{1}\in\binom{A\bs B}{2},\ D_{2}\in\binom{B\bs A}{2}\right\}.$$
	Therefore, we have $\left|\md_{\mf}(t)\right|\leq 2s+36\leq 11s+9t+18$.

	\medskip
	\noindent \textbf{Case 2.} $\left|A\cap B\right|= t-1$.
	\medskip

	Choose  $D\in\binom{[n]}{t+2}$ with $\left|D\cap A\right|\geq t$ and $\left|D\cap B\right|\geq t$. By (\ref{2509247}), we have $$\left|D\cap A\cap B\right|\geq \left|D\cap A\right|+\left|D\cap B\right|-(t+2)\geq t-2.$$ For convenience, write
		$$\md=\left\{ D_{1}\cup D_{2}\cup D_{3}: D_{1}\in\binom{A\cap B}{t-2},\ D_{2}\in\binom{A\bs B}{2},\ D_{3}\in\binom{B\bs A}{2}\right\}.$$
	 If $\left|D\cap A\cap B \right|=t-2$, then $D\in\md$.
	If $\left|D\cap A\cap B \right|=t-1$, then $D$ contains some members in 
	$\mw(A, B)$. For  $a\in A\bs B$, $b\in B\bs A$, write
	$\md(a,b)=\left\{D\in\md_{\mf}(t): (A\cap B)\cup\left\{a,b\right\}\subseteq D\right\}$. Therefore, we have
	 $$\md_{\mf}(t)\bs \left( \md_{\mf}(A;t)\cup \md_{\mf}(B;t)\right)\subseteq  \md\cup\left(\bigcup_{a\in A\bs B,\ b\in B\bs A}\md(a,b)\right).$$ 
	 
	 We claim $\left|\md(a,b)\right|\leq s+3$ for any $a\in A\bs B$ and $b\in B\bs A$. We may assume that $\md(a,b)\neq \emptyset$. Pick $D\in \md(a,b)$. Since $D\in \md_{\mf}(t)$ and $\min\left\{ \left|F_{1}\cap F_{2}\right|: F_{1}, F_{2}\in\md_{\mf}(t)\right\}=t-1$, there exists  $F\in\md_{\mf}(t)$ such that $\left|F\cap D\right|= t-1$. Hence $$\left|F\cap \left((A\cap B)\cup\left\{a,b\right\}\right)\right|\in\left\{t-2,t-1\right\}.$$
	When $\left|F\cap \left((A\cap B)\cup\left\{a,b\right\}\right)\right|=t-2$, we derive $\md(a,b)\subseteq \md_{\mf}(F;t)$, which implies $\left|\md(a,b)\right|\leq s$. When $\left|F\cap \left((A\cap B)\cup\left\{a,b\right\}\right)\right|=t-1$, there are at most three members in $\md(a,b)$ which are not in  $\md_{\mf}(F;t)$, implying that $\left|\md(a,b)\right|\leq s+3$.
	
	By the discussion above, we have $\left|\md_{\mf}(t)\right|\leq 2s+9(t-1)+9(s+3)=11s+9t+18$.
\end{proof}

\section{Proof of Theorem \ref{2509211}}\label{2510203}

In this section, we provide the proof of Theorem \ref{2509211}. Subsequently, in Remark \ref{2510291}, we demonstrate that the condition $n\geq L_{1}(k,t,s)$ of this theorem is optimal up to a constant factor.

\begin{proof} [\textnormal{\textbf{Proof of Theorem \ref{2509211}}}]
  Since  $\left\{F\in\binom{[n]}{k}: [t]\subseteq F\right\}$ is  $s$-almost $t$-intersecting, we have
\begin{equation}\label{2509217}
	\left|\mf\right|\geq \binom{n-t}{k-t}. 
\end{equation} 

\begin{cl}\label{2010161}
	$\mf$ is $t$-intersecting. 
\end{cl}
\begin{proof}
	Suppose for contradiction that $\mf$ is not $t$-intersecting. If $k=t+1$, then by Lemma \ref{2509215} (i) and Lemma \ref{2510072} (i), we have $\left|\mf\right|<n-t$, a contradiction to (\ref{2509217}). In the following, we assume $k\geq t+2$. 

	From Lemma \ref{2509214}, we obtain $\left|\md_{\mf}(t)\right|\leq s\binom{2k-2t+2}{k-t+1}$. Next, we show $$\left|\mi_{\mf}(t)\right|\leq \max\left\{g_{1}(n,k,t),g_{2}(n,k,t)\right\}.$$ If $| \bigcap_{F\in\mi_{\mf}(t)}F|<t$, then by Theorem \ref{2510151}, the  desired result holds. If $| \bigcap_{F\in\mi_{\mf}(t)}F|\geq t$, i.e., $\mi_{\mf}(t)\subseteq \{F\in\binom{[n]}{k}: X\subseteq F\}$ for some $X\in\binom{[n]}{t}$, then since  $\mf$ is not $t$-intersecting, we know $X\nsubseteq G$ for some $G\in\mf$. By $\left|G\bs X\right|\geq k-t+1$, there exists $H\in \binom{X\cup G}{k}$ such that $X\cap G\subseteq H$ and $\left|H\cap X\right|=t-1$. 
	For each $F\in\mi_{\mf}(t)$, we have  
	\begin{equation*}
		\begin{aligned}
			\left|F\cap \left(X\cup H\right)\right|&=\left|F\cap \left(X\cup G\right)\cap \left(X\cup H\right)\right|
			\geq  \left|F\cap \left(X\cup G\right) \right|+\left|X\cup H\right|-\left|X\cup G\right|\\
			&=\left|X\right|+\left|F\cap G\right|-\left|X\cap G\right|+\left|X\cup H\right|-\left|X\cup G\right|\\
			&=\left|F\cap G\right|+\left|X\cup H\right|-\left|G\right|
			\geq t+1.
		\end{aligned}
	\end{equation*}
	This implies 
	$$\left|\mi_{\mf}(t)\right|\leq \left| \left\{ F\in\binom{[n]}{k}: X\subseteq F,\ \left|F\cap \left( X\cup H\right)\right|\geq t+1\right\}\right|= g_{2}(n,k,t)-t.$$
	
	Therefore, we have $\left|\mf\right|\leq \max\left\{g_{1}(n,k,t),g_{2}(n,k,t)\right\}+s\binom{2k-2t+2}{k-t+1}$ when $k\geq t+2$. This together with Lemma \ref{2510072} (ii) yields $\left|\mf\right|<\binom{n-t}{k-t}$, a contradiction to (\ref{2509217}). 
\end{proof}

From (\ref{2509217}) and Lemma \ref{2510072}, we obtain $\left|\mf\right|>\max\left\{g_{1}(n,k,t),g_{2}(n,k,t)\right\}$. This combining with Theorem \ref{2510151} yields  $\mf \subseteq \left\{F\in\binom{[n]}{k}: X\subseteq F\right\}$ for some $t$-subset $X$ of $[n]$. By applying (\ref{2509217}) again, we have $\mf=\left\{F\in\binom{[n]}{k}: X\subseteq F\right\}$. \end{proof}

\begin{re}\label{2510291}
Let $n_{1}(k,t,s)$ be the minimum integer such that for any $n\geq n_{1}(k,t,s)$,  each $s$-almost $t$-intersecting family $\mf\subseteq\binom{[n]}{k}$ with the maximum size  consists of all $k$-subsets  containing a fixed $t$-subset of $[n]$. In the following, we show 
\begin{equation*}
	\frac{L_{1}(k,t,s)}{12e}\leq n_{1}(k,t,s)\leq L_{1}(k,t,s),
\end{equation*}
 where $e$ is  Euler’s number. Indeed, Theorem \ref{2509211} provides $n_{1}(k,t,s)\leq L_{1}(k,t,s)$.  

Let $n=(t+1)(k-t+1)$. Note that $\left\{F\in\binom{[n]}{k}: \left|F\cap [t+2]\right|\geq t+1\right\}$ is an $s$-almost $t$-intersecting family. This family is different from 
$\left\{F\in\binom{[n]}{k}: T\subseteq F\right\}$ for any $t$-subset $T$ of $[n]$, and its size is $\binom{n-t}{k-t}$. Hence $n_{1}(k,t,s)> (t+1)(k-t+1)$. 
 
Let $n=\lfloor \frac{(k-t+1)s^{\frac{1}{k-t}}}{2e} \rfloor$. For each positive integer $\ell$, by Stirling's formula, there exits a real number $\theta(\ell)\in (0,1)$ such that $\ell !=\sqrt{2\pi \ell}\left(\frac{\ell}{e}\right)^{\ell}e^{\frac{\theta(\ell)}{12\ell}}$, which implies $\ell !\geq \left(\frac{\ell}{e}\right)^{\ell}$. Then 
$$\binom{n-k-1}{k-t}\leq \frac{n^{k-t}}{(k-t)!}\leq \left(\frac{en}{k-t}\right)^{k-t}\leq \left( \frac{(k-t+1)s^{\frac{1}{k-t}}}{2(k-t)}\right)^{k-t}\leq s.$$
It follows that $\left\{F\in\binom{[n]}{k}: [t]\subseteq F\right\}\cup\left\{ [k+1]\bs\left\{1\right\}\right\}$ is $s$-almost $t$-intersecting. Therefore, we have  $n_{1}(k,t,s)> \frac{(k-t+1)s^{\frac{1}{k-t}}}{2e}-1$.

In summary, we know 
$n_{1}(k,t,s)\geq \frac{(t+1)(k-t+1)}{2}+\frac{(k-t+1)s^{\frac{1}{k-t}}}{4e}-\frac{1}{2}\geq \frac{L_{1}(k,t,s)}{12e}$, as desired.
\end{re}

\section{Proof of Theorem \ref{2509212}}\label{2510204}

An important tool for proving Theorem \ref{2509212} is the $t$-covers of families. Let $\mf\subseteq\binom{[n]}{k}$. 
A subset $T$  of $[n]$ is called a \textit{$t$-cover} of $\mf$ if $\left|T\cap F\right|\geq t$ for any $F\in \mf$. The \textit{$t$-covering number} $\tau_{t}(\mf)$ of $\mf$ is the minimum size of a $t$-cover of $\mf$. Suppose $\tau_{t}(\mf)\leq k$. Write 
$$\mt_{\mf}(t)=\left\{T\subseteq [n]: \left|T\right|\leq k,\ \textnormal{$T$ is a minimal (for containment) $t$-cover of $\mf$} \right\}.$$
For each positive $i$ with $\tau_{t}(\mf)\leq i\leq k$, denote
$$\mt_{\mf}(t,i)=\left\{T\in\mt_{\mf}(t): \left|T\right|=i\right\},\ \ \mf(i)=\left\{ F\in\binom{[n]}{k}: T\subseteq F\ \textnormal{for some $T\in \mt_{\mf}(t,i)$}\right\}.$$

The $s$-almost $t$-intersecting family $\mf$ is said to be \textit{maximal} if $\mg=\mf$ for any $s$-almost $t$-intersecting family $\mg$ with $\mf \subseteq \mg$. We will show a property about $t$-covers of such families. 

\begin{pr1}\label{2509221}
	Let $n$, $k$, $t$ and $s$ be positive integers with $k\geq t+1$ and $n\geq 2k-t+1$. Suppose that $\mf\subseteq \binom{[n]}{k}$ is a maximal $s$-almost $t$-intersecting family. If $\tau_{t}(\mf)\leq k$, then the following hold.
	\begin{enumerate}[\normalfont(i)]
		\item $\mt_{\mf}(t)$ is $t$-intersecting.
		\item $\mi_{\mf}(t)=\cup_{i=\tau_{t}(\mf)}^{k}\mf(i)$.
	\end{enumerate}
\end{pr1}
\begin{proof}
	(i) Suppose for contradiction that $\left|T_{1}\cap T_{2}\right|\leq t-1$ for some $T_{1},T_{2}\in\mt_{\mf}(t)$. By $\left|T_{2}\bs T_{1}\right|\geq \left|T_{2}\right|-t+1$, there exists $X\subseteq T_{2}\bs T_{1}$ such that $\left|X\right|=\left|T_{2}\right|-t+1$. Note that 
	$n-\left|T_{1}\cup X\right|\geq k-\left|T_{1}\right|$. Pick $Y\subseteq [n]\bs \left( T_{1}\cup X\right)$ with $\left|Y\right|=k-\left|T_{1}\right|$. Since $\mf$ is maximal and $T_{1}\subseteq T_{1}\cup Y$, we have $T_{1}\cup Y\in\mf$. But 
	$$\left|\left( T_{1}\cup Y\right)\cap T_{2}\right|=\left|T_{1}\cap T_{2}\right|+\left|Y\cap T_{2}\right|\leq \left|T_{1}\cap T_{2}\right|+\left| \left(T_{2}\bs T_{1}\right)\bs X\right|=t-1.$$
	This contradicts to $T_{2}\in\mt_{\mf}(t)$.
  Then $\mt_{\mf}(t)$ is $t$-intersecting.
	
	(ii) For each $F\in$ LHS, we know $F$ is a $t$-cover of $\mf$. Then $F$ contains a minimal $t$-cover of $\mf$ whose size is no more than $k$. Hence LHS $\subseteq$ RHS. For each $F\in$ RHS, since $\mf$ is maximal, we have $F\in\mf$. Moreover, by the definition of $\mi_{\mf}(t)$, we know $F\in$ LHS. Hence RHS $\subseteq$ LHS. The desired result holds.
\end{proof}

For a considered family $\mf$, by Proposition \ref{2509221} (ii), we can  estimate the size of $\mi_{\mf}(t)$ in terms of the sizes of families $\mf(\tau_{t}(\mf)), \mf(\tau_{t}(\mf)+1), \ldots,\mf(k)$.  Moreover, the size of  $\mf(i)$ can be bounded by the size of $\mt_{\mf}(t,i)$.  The subsequent two results  are related to the sizes of $\mt_{\mf}(t,i)$ and $\mf(i)$, respectively.

\begin{lem}\label{2509222}
	Let $n$, $k$, $t$ and $s$ be positive integers with $k\geq t+1$ and $n\geq 2k-t+1$.  Suppose that $\mf\subseteq \binom{[n]}{k}$ is a maximal $s$-almost $t$-intersecting family. If $\tau_{t}(\mf)\leq k$, then 
	$$\sum_{i=\tau_{t}(\mf)}^{k}(k-2t+\tau_{t}(\mf))^{k-i}\left|\mt_{\mf}(t,i)\right|\leq (k-2t+\tau_{t}(\mf))^{k-t}\binom{\tau_{t}(\mf)}{t}.$$
\end{lem}

\begin{proof}
	For a sequence $S=(x_{1},x_{2},\ldots,x_{i})$, where $x_{1},x_{2},\ldots,x_{i}$ are distinct,  denote $\widehat{S}=\left\{x_{1},x_{2},\ldots,x_{i} \right\}$. 
	 Set $m=k-2t+\tau_{t}(\mf)$, and $$Y_{i+1}=\left\{ -(i-t+1)m+1, -(i-t+1)m+2,\ldots, -(i-t)m\right\},\ t\leq i\leq k-1.$$

	Now, we construct a series of sequences. First, fix $T_{0}\in\mt_{\mf}(t,\tau_{t}(\mf))$. For each $t$-subset of $T_{0}$, we fix an ordering of its elements and define a sequence accordingly. The family consisting of such $\binom{\tau_{t}(\mf)}{t}$ sequences is denoted by $\ms(t)$.  
	
	Next, we build step by step the family $\ms(i+1)$ consisting of sequences with length $i+1$, where $t\leq i\leq k-1$. Suppose that $\ms(i)$  have been built. Pick $S=\left(x_{1},x_{2},\ldots,x_{i}\right)\in\ms(i)$.  We set $\my_{1}(S)$ as following. If $\mid\widehat{S}\cap F \mid\geq t$ for any $F\in\mf$, then set
	$$\my_{1}(S)=\left\{(x_{1},x_{2},\ldots,x_{i},y): y\in Y_{i+1}\right\}.$$  Otherwise, we pick $F\in\mf$ with $\mid\widehat{S}\cap F \mid<t$ and set  
	$$\my_{1}(S)=\left\{(x_{1},x_{2},\ldots,x_{i},y): y\in F\bs \widehat{S}\right\}.$$ 
	Define $\ms(i+1)=\bigcup_{S\in\ms(i)}\my_{1}(S)$. 
	 
	For each $S\in\cup_{i=t}^{k-1}\ms(i)$, there exists a $t$-subset $Z$ of $T_{0}$ such that $Z\subseteq \widehat{S}$,  which implies for each $F\in\mf$
	$$\mid F\bs\widehat{S}\mid\leq \left|F\right|-\mid F\cap Z\mid=k-t+\left|\left(F\cap T_{0}\right)\cup Z\right|-\left|F\cap T_{0}\right|\leq k-2t+\tau_{t}(\mf).$$
	Therefore, we have $\left|\ms(i+1)\right|\leq m\left|\ms(i)\right|$ for any $t\leq i\leq k-1$. We further derive 
	\begin{equation}\label{2509231}
		\left|\ms(k)\right|\leq m ^{k-t}\left|\ms(t)\right|=m^{k-t}\binom{\tau_{t}(\mf)}{t}.
	\end{equation}
	
	\begin{cl}
	 For each  $T\in\mt_{\mf}(t)$, there exists $S\in\ms(\left|T\right|)$ such that $\widehat{S}=T$.
	\end{cl}
	\begin{proof}
		By Lemma \ref{2509221} (i), we have $\left|T\cap T_{0}\right|\geq t$. Then $\{S\in\ms(t): \widehat{S}\subseteq T\}$ is non-empty. Let $i$ be the maximum positive integer such that $\{S\in\ms(i): \widehat{S}\subseteq T\}\neq \emptyset$, which implies $i\leq \left|T\right|$. It is sufficient to show $i=\left|T\right|$.
		
		Suppose for contradiction that $i<\left|T\right|$. Pick $S=(x_{1},x_{2},\ldots,x_{i})\in\ms(i)$ with $ \widehat{S}\subsetneq T$.  By the minimality of $T$, we know $\widehat{S}$ is not a $t$-cover of $\mf$. By the definition of $\ms(i+1)$, there exists $F\in \mf$ such that $\mid\widehat{S}\cap F\mid<t$ and 
		$$\left\{(x_{1},x_{2},\ldots,x_{i},y): y\in F\bs \widehat{S}\right\}\subseteq \ms(i+1).$$
		From $\left|T\cap F\right|\geq t$, we know $T\cap (F\bs \widehat{S})\neq \emptyset$. Therefore, we obtain $\{S\in\ms(i+1): \widehat{S}\subseteq T\}\neq \emptyset$, a contradiction to the choice of $i$. Then $i=\left|T\right|$, the desired result holds.
	\end{proof}
	
	Pick $T\in\mt_{\mf}(t)$. By the claim above, there exists $S=(x_{1},x_{2},\ldots,x_{\left|T\right|})\in\ms(\left|T\right|)$ such that $\widehat{S}=T$. Set 
	$$\my_{2}(T)=\left\{(x_{1},x_{2},\ldots,x_{\left|T\right|},y_{\left|T\right|+1},\ldots,y_{k}):y_{\left|T\right|+1}\in Y_{\left|T\right|+1},\ldots, y_{k}\in Y_{k}\right\}.$$
	From the definition of $\ms(k)$, we obtain $\my_{2}(T)\subseteq \ms(k)$.
	Moreover, if $T_{1}$ and $T_{2}$ are distinct members in $\mt_{\mf}(t)$, then $\my_{2}(T_{1})\cap \my_{2}(T_{2})=\emptyset$. Otherwise, we have $T_{1}\subsetneq T_{2}$ or $T_{1} \supsetneq T_{2}$, a contradiction.

	Therefore, we have 
	$\sum_{i=\tau_{t}(\mf)}^{k}m^{k-i}\left|\mt_{\mf}(t,i)\right|=|\bigcup_{T\in \mt_{\mf}(t)}\my_{2}(T)|\leq \left|\ms(k)\right|$.
	This together with (\ref{2509231}) yields the required result.
\end{proof}

Some frequently used functions are presented below:
\begin{equation}\label{2510082}
	f_{1}(n,k,t,x,y)=(k-2t+x)^{y-t}\binom{x}{t}\binom{n-y}{k-y},
\end{equation}
\begin{equation}\label{2509295}
	f_{2}(n,k,t,s,x)=\binom{n-t}{k-t}-\binom{n-t-x}{k-t}+\left(k-t-x+1\right)^{2}\binom{n-t-2}{k-t-2}+2s, 
\end{equation}
\begin{equation}\label{2509296}
	f_{3}(n,k,t,s)=f_{2}(n,k,t,s,k-t)+t(k-t)\binom{n-t-2}{k-t-2}+s.
\end{equation}

\begin{cor}\label{2509242}
Let $n$, $k$, $t$ and $s$ be positive integers with $k\geq t+1$ and $n\geq (t+1)(k-t+1)^{2}$. Suppose that $\mf\subseteq \binom{[n]}{k}$ is a maximal $s$-almost $t$-intersecting family. If $\tau_{t}(\mf)\leq k$, then for each  $j\in\left\{\tau_{t}(\mf),\tau_{t}(\mf)+1,\ldots,k\right\}$, we have  $\left|\cup_{i=j}^{k}\mf(i)\right|\leq f_{1}(n,k,t,\tau_{t}(\mf), j)$.
\end{cor}
\begin{proof}	
	Set $m=k-2t+\tau_{t}(\mf)$. Then
	\begin{equation}\label{2511181}
		\left|\cup_{i=j}^{k}\mf(i)\right|\leq \sum_{i=j}^{k}\left|\mt_{\mf}(t,i)\right|\binom{n-i}{k-i}=\sum_{i=j}^{k}m^{k-i}\cdot \frac{\binom{n-i}{k-i}}{m^{k-i}}\cdot\left|\mt_{\mf}(t,i)\right|.
	\end{equation}

	We  claim that $\frac{1}{m^{k-j}}\binom{n-j}{k-j}\geq \frac{1}{m^{k-i}}\binom{n-i}{k-i}$ for any $i\in \left\{ j,j+1,\ldots,k\right\}$. If $j=k$, then there is nothing to prove. If $j\leq k-1$, then for each $i\in \left\{j,j+1,\ldots,k-1\right\}$
	$$\frac{\binom{n-i}{k-i}}{m^{k-i}}\cdot \frac{m^{k-i-1}}{\binom{n-i-1}{k-i-1}}= \frac{1}{m}\cdot \frac{n-i}{k-i}\geq \frac{1}{2(k-t)}\cdot \frac{ (t+1)(k-t+1)^{2}-t}{k-t}\geq 1, $$
	which implies the required claim.
	It follows from Lemma \ref{2509222} and (\ref{2511181}) that
	$$\left|\cup_{i=j}^{k}\mf(i)\right|\leq \frac{1}{m^{k-j}}\binom{n-j}{k-j}\sum_{i=j}^{k}m^{k-i}\left|\mt_{\mf}(t,i)\right|\leq m^{j-t} \binom{\tau_{t}(\mf)}{t}\binom{n-j}{k-j}, $$
	as desired.
\end{proof}

In the following, we  determine the maximum-sized $s$-almost $t$-intersecting family, subject to the condition that it is not $t$-intersecting. We first show an example.

\begin{ex}\label{2509243}
	Let $n$, $k$, $t$ and $s$ be positive integers with $\binom{n-k-1}{k-t}\geq s$. Suppose  $X\in\binom{[n]}{k+1}$ and $Y\in\binom{X}{t}$.  For any $s$-subset $\ma$ of $\left\{ F\in\binom{[n]}{k}: F\cap X=Y\right\}$ and $\min\left\{t,s\right\}$-subset $\mb$ of $\left\{ F\in\binom{X}{k}: Y \nsubseteq F\right\}$, the family $$\mh(X,Y;k,t)\cup\ma\cup\mb$$  is $s$-almost $t$-intersecting but not $t$-intersecting. 
\end{ex}

The size of the family in Example \ref{2509243} is 
\begin{equation}\label{2509244}
	\binom{n-t}{k-t}-\binom{n-k-1}{k-t}+s+\min\left\{t,s\right\}=:g_{3}(n,k,t,s).
\end{equation}
This can be used to bound the sizes of  some families that have a large $t$-covering number.

\begin{lem}\label{2509245}
	Let $n,k,t$ and $s$ be positive integers with $k\geq t+2$ and $n\geq L_{2}(k,t,s)$. Suppose that $\mf\subseteq \binom{[n]}{k}$ is a maximal $s$-almost $t$-intersecting family but not $t$-intersecting. If $\tau_{t}(\mf)\geq t+2$, then $\left|\mf\right|< g_{3}(n,k,t,s)$.
\end{lem}
\begin{proof}
	We first show 
	\begin{equation}\label{2510113}
		\left|\mf\right|\leq f_{1}(n,k,t,t+2,t+2) +\frac{11}{54}(k-t+1)\binom{n-t-1}{k-t-1}.
	\end{equation}

	By Lemma \ref{2509214} and Lemma \ref{2509215} (ii), we have $\left|\md_{\mf}(t)\right|\leq 11s+9t+18$ when $k=t+2$; $\left|\md_{\mf}(t)\right|\leq s \binom{2k-2t+2}{k-t+1}$ when $k\geq t+3$. This together with Lemma \ref{2510071} (ii) yields $\left|\md_{\mf}(t)\right|\leq \frac{11}{54}(k-t+1)\binom{n-t-1}{k-t-1}$. 	
	If  $\tau_{t}(\mf)\geq k+1$, then $\mf=\md_{\mf}(t)$. Hence (\ref{2510113}) holds. If $\tau_{t}(\mf)\leq k$, then by Proposition \ref{2509221} (ii), we have $\mi_{\mf}(t)=\cup_{i=\tau_{t}(\mf)}^{k}\mf(i)$. It follows from Corollary \ref{2509242} that $\left|\mi_{\mf}(t)\right|\leq f_{1}(n,k,t,\tau_{t}(\mf),\tau_{t}(\mf))$. From Lemma \ref{2510081}, we get $\left|\mi_{\mf}(t)\right|\leq f_{1}(n,k,t,t+2,t+2)$. This combining with $\left|\mf\right|=\left|\mi_{\mf}(t)\right|+\left|\md_{\mf}(t)\right|$ implies (\ref{2510113}). 
	
	According to (\ref{2510113}), Lemma \ref{2510086} (i), we get $\left|\mf\right|< g_{3}(n,k,t,s)$.
\end{proof}

Now, we are ready to prove Theorem \ref{2509212}.

\begin{proof}[\textnormal{\textbf{Proof of Theorem \ref{2509212}}}]

Note that $L_{2}(k,t,s)\geq L_{1}(k,t,s)$. From Lemma \ref{2510071} (i), we obtain $\binom{n-k-1}{k-t}\geq s$. Since $\mf$ has the maximum size, by Example \ref{2509243}, we have
\begin{equation}\label{2509292}
	\left|\mf\right|\geq  g_{3}(n,k,t,s).
\end{equation}
Moreover, we know that $\mf$ is maximal. According to Lemma \ref{2509245}, we get $\tau_{t}(\mf)\leq t+1$. Note that $\mf$ is not $t$-intersecting. Hence $\tau_{t}(\mf)=t+1$.

\begin{cl}\label{2509291}
	The following hold.
	\begin{enumerate}[\normalfont(i)]
		\item 	$\left|\mt_{\mf}(t,t+1)\right|\geq 2$.
		\item   There exists a $t$-subset $Y$ of $[n]$ such that $Y\subseteq T$ for any $T\in\mt_{\mf}(t,t+1)$. 
	\end{enumerate}
\end{cl}
\begin{proof} (i)
	Suppose for contradiction that $\mt_{\mf}(t,t+1)=\left\{T\right\}$.  From Proposition \ref{2509221} (ii) and Corollary \ref{2509242}, we obtain
	$$\left|\mi_{\mf}(t)\right|\leq\left|\left\{F\in\mf: T\subseteq F\right\}\right|+\left|\cup_{i=t+2}^{k}\mf(i)\right|\leq \binom{n-t-1}{k-t-1}+f_{1}(n,k,t,t+1,t+2).$$
	By Lemma \ref{2509214}, Lemma \ref{2509215} (ii) and Lemma \ref{2510071} (ii), we get $\left|\md_{\mf}(t)\right|<\frac{11}{54}(k-t+1)\binom{n-t-1}{k-t-1}$. Applying Lemma \ref{2510086} (ii), we know $\left|\mf\right|<  g_{3}(n,k,t,s)$ a contradiction to (\ref{2509292}). 
	
	(ii) By (i), there exist two distinct  $T_{1}, T_{2}$ in $\mt_{\mf}(t,t+1)$. This together with Proposition \ref{2509221} (i) yields $\left|T_{1}\cap T_{2}\right|=t$. Next, we show $Y:=T_{1}\cap T_{2}\subseteq T$ for any $T\in \mt_{\mf}(t,t+1)$.
	
	Suppose for contradiction that $Y\nsubseteq T$ for some $T\in\mt_{\mf}(t,t+1)$. Set $T_{1}\bs T_{2}=\left\{x_{1}\right\}$ and $T_{2}\bs T_{1}=\left\{x_{2}\right\}$. By Proposition \ref{2509221} (i), we have $\left|T\cap T_{1}\right|\geq t$ and $\left|T\cap T_{2}\right|\geq t$, implying that $\left|T\cap Y\right|=t-1$ and $T=\left(T\cap Y\right)\cup\left\{x_{1},x_{2}\right\}$. Recall that $\mf$ is not $t$-intersecting. It follows from $\tau_{t}(\mf)=t+1$ that $\left|F_{1}\cap F_{2}\right|=t-1$ for some $F_{1},F_{2}\in\mf$. Then for  $i\in\left\{1,2\right\}$
	$$t-1\geq \left|T_{i}\cap F_{1}\cap F_{2}\right|= \left|T_{i}\cap F_{1}\right|+\left|T_{i}\cap F_{2}\right|-\left|T_{i}\cap \left(F_{1}\cup F_{2}\right)\right|\geq t-1.$$
	This implies $F_{1}\cap F_{2}\subseteq Y$. Set $Y\bs\left(F_{1}\cap F_{2}\right)=\left\{x_{0}\right\}$.
	
	Note that $\left\{x_{0}, x_{1}, x_{2}\right\}\cap \left(F_{1}\cap F_{2}\right)=\emptyset$. By $\left|T_{1}\cap F_{1} \right|\geq t$ and $\left|T_{1}\cap F_{2}\right|\geq t$, we have 
	$$\left| \left(F_{1}\cap F_{2}\right)\cup \left( \left\{x_{0},x_{1}\right\}\cap \left(F_{1}\bs F_{2}\right)\right)\right|=\left|T_{1}\cap F_{1}\right|\geq t,$$
	$$\left| \left(F_{1}\cap F_{2}\right)\cup \left( \left\{x_{0},x_{1}\right\}\cap \left(F_{2}\bs F_{1}\right)\right)\right|=\left|T_{1}\cap F_{2}\right|\geq t.$$
	Therefore $x_{0}\in F_{1}\bs F_{2}$ and $x_{1}\in F_{2}\bs F_{1}$, or $x_{0}\in F_{2}\bs F_{1}$ and $x_{1}\in F_{1}\bs F_{2}$.  W.l.o.g.,  assume that $x_{0}\in F_{1}\bs F_{2}$ and $x_{1}\in F_{2}\bs F_{1}$. This together with
	$$\left| \left(F_{1}\cap F_{2}\right)\cup \left( \left\{x_{0},x_{2}\right\}\cap \left(F_{2}\bs F_{1}\right)\right)\right|=\left|T_{2}\cap F_{2}\right|\geq t$$
	 implies $x_{2}\in F_{2}\bs F_{1}$. Then 
	$$\left|T\cap F_{1}\right|=\left|\left(T\cap Y\cap F_{1}\right)\cup \left(\left\{x_{1},x_{2} \right\}\cap F_{1}\right) \right|=\left|T\cap Y\right|=t-1,$$
	a contradiction to $T\in \mt_{\mf}(t,t+1)$. Then $Y\subseteq T$ for any $T\in \mt_{\mf}(t,t+1)$.
\end{proof}

Set $X=\bigcup_{T\in \mt_{\mf}(t,t+1)} T$. Then $Y\subseteq X$.  By Claim \ref{2509291}, we have $\left|X\bs Y\right|\geq 2$ and $\mt_{\mf}(t,t+1)=\left\{Y\cup\left\{x\right\}: x\in X\bs Y\right\}$. In the remaining of this proof, write 
$$\mf_{1}=\left\{F\in\mf: Y\subseteq F\right\},\ \mf_{2}=\left\{F\in \mf: Y\nsubseteq F\right\}.$$
Since $\mf$ is not $t$-intersecting, we know $\mf_{2}\neq \emptyset$. Moreover, we have 
\begin{equation}\label{2509293}
	\left|F\cap Y\right|=t-1\ \textnormal{and}\ X\bs Y\subseteq F\ \textnormal{for any $F\in\mf_{2}$}.
\end{equation}
It follows that $\mf_{2}$ is $t$-intersecting.  Then there exist $A \in \mf_{1}$ and $B\in \mf_{2}$ such that $\left|A\cap B\right|<t$. Indeed, we know $\left| A\cap B\right|=t-1$ and $A\cap \left( X\bs Y\right)=\emptyset$ due to (\ref{2509293}).

\begin{cl}\label{2509294}
	$F\bs Y=B\bs Y$ for any $F\in \mf_{2}$.
\end{cl}
\begin{proof}
	Suppose for contradiction that $C\bs Y\neq B\bs Y$ for some $C\in \mf_{2}$. Set $B^{\prime}=B\bs Y$ and $C^{\prime}=C\bs Y$. By (\ref{2509293}), we know $\left|B^{\prime}\right|=\left|C^{\prime}\right|=k-t+1$.
	
	Pick $F\in\mf_{1}$. If $\left|F\cap B\right|\geq t$ and $\left|F\cap C\right|\geq t$, then $F\cap B^{\prime}$ and $F\cap C^{\prime}$ are non-empty, implying that $F\cap B^{\prime}\cap C^{\prime}\neq \emptyset$, or $F\cap \left(B^{\prime}\bs C^{\prime}\right)\neq \emptyset$ and $F\cap \left( C^{\prime}\bs B^{\prime}\right)\neq \emptyset$. Hence
	\begin{equation*}
		\begin{aligned}
			&\left|\mf_{1}\bs \left(\md_{\mf}(B;t)\cup \md_{\mf}(C;t)\right)\right|\\
			\leq&  \binom{n-t}{k-t}-\binom{n-t-\left|B^{\prime}\cap C^{\prime}\right|}{k-t}+\left|B^{\prime}\bs C^{\prime}\right|\left|C^{\prime}\bs B^{\prime}\right|\binom{n-t-2}{k-t-2}. 
		\end{aligned}
	\end{equation*}
	Note that $\left|B^{\prime}\bs C^{\prime}\right|=\left|C^{\prime}\bs B^{\prime}\right|=k-t-\left|B^{\prime}\cap C^{\prime}\right|+1$.
	Then $\left|\mf_{1}\right|\leq f_{2}(n,k,t,s, \left|B^{\prime}\cap C^{\prime}\right|)$. 
	By $B^{\prime}\neq C^{\prime}$, we have $\left|B^{\prime}\cap C^{\prime}\right|\leq k-t$. It follows from Lemma \ref{2510114} that $\left|\mf_{1}\right|\leq f_{2}(n,k,t,s, k-t)$.	
	
	Choose $F\in\mf_{2}$. If $\left|F\cap A\right|\geq t$, then $F\cap \left(A\bs Y\right)\neq \emptyset$. 	This together with (\ref{2509293}) implies 
	$$\mf_{2}\bs \md_{\mf}(A;t)\subseteq \left\{F\in\binom{[n]}{k}: \left| F\cap Y\right|=t-1,\ X\bs Y\subseteq F,\ F\cap \left(A\bs Y\right)\neq\emptyset \right\}.$$
	Recall that $\left|X\bs Y\right|\geq 2$ and $A\cap \left( X\bs Y\right)=\emptyset$. We derive 
	$$\left|\mf_{2}\bs \md_{\mf}(A;t)\right|\leq t(k-t)\binom{n-t-\left|X\bs Y\right|}{k-t-\left|X\bs Y\right|}\leq t(k-t)\binom{n-t-2}{k-t-2}. $$
	It follows from $\mf=\mf_{1}\cup \mf_{2}$ that $\left|\mf\right|\leq f_{3}(n,k,t,s)$. By Lemma \ref{2510086} (iii), we have $f_{3}(n,k,t,s)< g_{3}(n,k,t,s)$, a contradiction to 
	(\ref{2509292}).
\end{proof}
From (\ref{2509293}), we obtain $\left|B\cap Y\right|=t-1$. Then 
$$\left\{F\in \mf_{1}: \left|F\cap B\right|\geq t\right\}\subseteq \mh(Y\cup B, Y;k,t),$$
$$\left\{F\in \mf_{1}: \left|F\cap B\right|< t\right\}\subseteq \left\{F\in\binom{[n]}{k}: F\cap\left(Y\cup B\right)=Y\right\}\cap \md_{\mf}(B;t).$$
For each $F\in \mf_{2}$, by Claim \ref{2509294}, we know $F\subseteq Y\cup \left(F\bs Y\right)=Y\cup B$. It follows from $A\cap \left(Y\cup B\right)=Y$ that 
$\left|F\cap A\right|=\left| F\cap\left(Y\cup B\right)\cap A \right|=\left|F\cap Y\right|<t$. Therefore 
$$\mf_{2}\subseteq \left\{F\in\binom{Y\cup B}{k}: Y\nsubseteq F\right\}\cap \md_{\mf}(A;t).$$
These together with (\ref{2509292}) yield
\begin{equation*}
	\begin{aligned}
		 g_{3}(n,k,t,s)=&\binom{n-t}{k-t}-\binom{n-k-1}{k-t}+s+\min\left\{t,s\right\}\\
		\geq& \left|\left\{F\in \mf_{1}: \left|F\cap B\right|\geq t\right\}\right|+s+\min\left\{t,s\right\}\\
		\geq& \left|\left\{F\in \mf_{1}: \left|F\cap B\right|\geq t\right\}\right|+\left|\left\{F\in \mf_{1}: \left|F\cap B\right|< t\right\}\right|+\min\left\{t,s\right\}\\
		\geq& \left|\left\{F\in \mf_{1}: \left|F\cap B\right|\geq t\right\}\right|+\left|\left\{F\in \mf_{1}: \left|F\cap B\right|< t\right\}\right|+\left|\mf_{2}\right|\geq  g_{3}(n,k,t,s),
	\end{aligned}
\end{equation*}
which implies $\left\{F\in \mf_{1}: \left|F\cap B\right|\geq t\right\}= \mh(Y\cup B, Y;k,t)$, $\left\{F\in \mf_{1}: \left|F\cap B\right|< t\right\}$ is an $s$-subset of  $\left\{F\in\binom{[n]}{k}: F\cap\left(Y\cup B\right)=Y\right\}$, $\mf_{2}$ is a $\min\left\{t,s\right\}$-subset of $\left\{F\in\binom{Y\cup B}{k}: Y\nsubseteq F\right\}$.
This finishes the  proof.
\end{proof}

\section{Proof of Theorem \ref{2509213}} \label{2510205}

This section is devoted to the proof of Theorem \ref{2509213}. It is worth noting that the condition $s\geq 7$ imposed in  this theorem is necessary, as detailed in Remark \ref{2601021}. Turning to the main argument, we begin  by presenting two examples.

\begin{ex}\label{2509301}
	Let $n$, $t$ and $s$ be positive integers with $n\geq t+s+2$. Suppose   $X\in \binom{[n]}{t+s+2}$, $Y\in\binom{X}{t+1}$ and  $Z\in\binom{Y}{t-1}$. The family $\mm_{1}(X,Y,Z;t)$  is $s$-almost $t$-intersecting but not $t$-intersecting.
\end{ex}

\begin{ex}\label{2509302}
	Let $n$, $t$ and $s$ be positive integers with $t\geq s$ and $n\geq t+3$. Suppose    $X\in \binom{[n]}{t+3}$, $Y\in\binom{X}{2}$ and $Z\in\binom{X\bs Y}{t-s}$.  The family $\mm_{2}(X,Y,Z;t)$ is $s$-almost $t$-intersecting but not $t$-intersecting.
\end{ex}

The sizes of the families in Examples \ref{2509301} and \ref{2509302}  are equal to $2s+3$. Next, we show some properties concerning  families whose sizes are at least $2s+3$.
Suppose $\mf\subseteq \binom{[n]}{t+1}$, $A, B\in \binom{[n]}{t+1}$ and $c\in [n]$. For convenience,  write 
$$\mi_{\mf}(A;t)=\left\{ F\in\mf: \left|F\cap A\right|=t\right\}, $$
$$\mmu(A,B;c)=\left\{\left(A\cap B\right)\cup \left\{c,x\right\}: x\in [n]\bs \left(A\cup B\right)\right\},$$
$$\mv(A,B;c)=\left\{\left(A\cup B\right)\bs \left\{c,x\right\}: x\in A\cap B\right\}.$$  
 In the following, we  also make use of the notation $\mw(A, B)$ introduced  in (\ref{2512291}).                                       
\begin{pr1}\label{2509303}
	Let $n$, $t$ and $s$ be positive integers with $n\geq t+3$. Suppose that $\mf\subseteq \binom{[n]}{t+1}$ is  $s$-almost $t$-intersecting  but not $t$-intersecting. If $\left|\mf\right|\geq 2s+3$, then the following hold.
	\begin{enumerate}[\normalfont(i)]
		\item $\left| \mi_{\mf}(F;t)\right|\geq s+2$ for any $F\in\mf$.
		\item $\left|A\cap B\right|\geq t-1$ for any $A, B\in\mf$.
		\item  $\left\{B\right\}\cup \left(\mi_{\mf}(B;t)\bs\mw(A,B)\right)\subseteq \md_{\mf}(A;t)$ for any $A, B\in\mf$ with $\left|A\cap B\right|=t-1$.
	\end{enumerate}
\end{pr1}
\begin{proof}
	(i) For each $F\in\mf$,  we know $\mf=\left\{F\right\}\cup \mi_{\mf}(F;t)\cup \md_{\mf}(F;t)$. It follows from $\left|\mf\right|\geq 2s+3$ and $\left|\md_{\mf}(F;t)\right|\leq s$ that $\left|\mi_{\mf}(F;t)\right|\geq s+2$.
	
	(ii) Suppose for contradiction that $\left|A\cap B\right|\leq t-2$ for some $A,B\in \mf$. Then $\mf=\md_{\mf}(A;t)\cup \md_{\mf}(B;t)$, which implies $\left|\mf\right|\leq 2s$, a contradiction.
	
	(iii)  Observe that  $\mi_{\mf}(B;t)\bs \md_{\mf}(A;t)\subseteq \mw(A,B)$ due to (\ref{2510021}). 
	This implies $\mi_{\mf}(B;t)\bs\mw(A,B)\subseteq \md_{\mf}(A;t)$. Then the desired result holds.
\end{proof}

\begin{pr1}\label{2510031}
	Let $n$, $t$ and $s$ be positive integers with $s\geq 7$ and $n\geq t+3$. Suppose that $\mf\subseteq \binom{[n]}{t+1}$ is $s$-almost $t$-intersecting  but not $t$-intersecting.   If $\left|\mf\right|\geq 2s+3$, then $\md_{\mf}(F;t)$ is $t$-intersecting for any $F\in \mf$ with $\md_{\mf}(F;t)\neq \emptyset$.
\end{pr1}
\begin{proof}
 Let $F\in \mf$ with $\md_{\mf}(F;t)\neq \emptyset$.	Suppose for contradiction that $\left|A\cap B\right|<t$ for some $A, B\in\md_{\mf}(F;t)$.  This together with Proposition \ref{2509303} (ii) yields $\left|F\cap A\right|=\left|F\cap B\right|=\left|A\cap B\right|=t-1$. It follows from Proposition \ref{2509303} (iii) that 
	$$\left\{A,B\right\}\cup \left( \left(\mi_{\mf}(A;t)\cup \mi_{\mf}(B;t) \right)\bs \left(\mw(F,A)\cup \mw(F,B) \right)  \right)\subseteq \md_{\mf}(F;t).$$
	Notice that $\mi_{\mf}(A;t)\cap \mi_{\mf}(B;t)\subseteq \mw(A,B)$ due to (\ref{2510021}). From Proposition \ref{2509303} (i), we obtain
	\begin{equation*}
		\begin{aligned}
			\left|\md_{\mf}(F;t)\right|&\geq 2+\left|\mi_{\mf}(A;t)\cup \mi_{\mf}(B;t)\right|-\left|\mw(F,A)\cup \mw(F,B)\right|\\
			&\geq  2+\left|\mi_{\mf}(A;t)\right|+\left|\mi_{\mf}(B;t)\right| -\left|\mw(A,B)\right|-\left|\mw(F,A)\right|-\left|\mw(F,B)\right|
			\geq s+1,
		\end{aligned}
	\end{equation*}
	a contradiction  to the fact that $\mf$ is $s$-almost $t$-intersecting.
\end{proof}

\begin{lem}\label{2510311}
	Let $n$, $t$ and $s$ be positive integers with $s\geq 7$ and $n\geq t+3$. Suppose that $\mf\subseteq \binom{[n]}{t+1}$ is $s$-almost $t$-intersecting  but not $t$-intersecting. The following hold.
	\begin{enumerate}[\normalfont(i)]
		\item If $\left|\mf\right|\geq 2s+3$ and $\left|\bigcup_{F\in\mf}F\right|\geq t+4$, then $n\geq t+s+2$ and $\mf=\mm_{1}(X,Y,Z;t)$ for some $X\in \binom{[n]}{t+s+2}$, $Y\in\binom{X}{t+1}$ and $Z\in\binom{Y}{t-1}$.
		\item If $\left|\mf\right|\geq 2s+3$ and $\left|\bigcup_{F\in\mf}F\right|= t+3$, then $t\geq s$ and $\mf=\mm_{2}(X,Y,Z;t)$ for some $X\in \binom{[n]}{t+3}$, $Y\in\binom{X}{2}$ and $Z\in\binom{X\bs Y}{t-s}$.
 	\end{enumerate}
\end{lem}

\begin{proof} 
	Assume that $\mf\subseteq \binom{[n]}{t+1}$ is an $s$-almost $t$-intersecting family that is  not $t$-intersecting, and $\left|\mf\right|\geq 2s+3$.
	Let $A, B\in \mf$ with $\left|A\cap B\right|<t$. By Proposition \ref{2509303} (ii), we have $\left|A\cap B\right|=t-1$.    In this proof, set  $A\bs B=\left\{a,c\right\}$ and $B\bs A=\left\{b,d\right\}$, and write
	$$\mf_{1}=\left\{B\right\}\cup \left\{ F\in\mf: \left|F\cap A\right|= t-1,\ \left|F\cap B\right|= t\right\}, $$
	$$\mf_{2}=\left\{A\right\}\cup \left\{ F\in\mf: \left|F\cap A\right|= t,\ \left|F\cap B\right|=t-1\right\}, $$
	$$\mf_{3}=\left\{F\in\mf: \left|F\cap A\right|\geq t,\ \left|F\cap B\right|\geq t \right\}.$$
   From Proposition \ref{2509303} (ii) and Proposition \ref{2510031}, we obtain  $\mf_{1}=\md_{\mf}(A;t)$ and $\mf_{2}=\md_{\mf}(B;t)$, which implies $\mf=\mf_{1}\cup \mf_{2}\cup \mf_{3}$. Note that $\mf_{3}\subseteq \mw(A, B)$ due to (\ref{2510021}). This together with $\left|\mf\right|\geq 2s+3$, $\left|\mf_{1}\right|\leq s$ and $\left|\mf_{2}\right|\leq s$ yields
	\begin{equation}\label{2510033}
		\left|\mf_{1}\bs \left\{B\right\}\right|=\left|\mf\right|-\left|\mf_{2}\right|-\left|\mf_{3}\right|-1\geq s-2,\ \left|\mf_{2}\bs\left\{ A\right\}\right|=\left|\mf\right|-\left|\mf_{1}\right|-\left|\mf_{3}\right|-1\geq s-2.
	\end{equation}

(i) Note that the members in $\mf_{3}\cup \left\{A, B\right\}$ are subsets of $A\cup B$. It follows from $\left|\bigcup_{F\in\mf}F\right|\geq t+4$ that  there exists $F\in \left(\mf_{1}\bs \left\{B\right\}\right)\cup \left( \mf_{2}\bs \left\{A\right\}\right)$ such that $F\nsubseteq A\cup B$. Moreover, we know  
	\begin{equation}\label{2510042}
		\textnormal{if}\ F\in \left(\mf_{1}\bs \left\{B\right\}\right)\cup \left( \mf_{2}\bs \left\{A\right\}\right)\ \textnormal{with}\   F\nsubseteq A\cup B,\ \textnormal{then}\ A\cap B\subseteq F.
	\end{equation}
	
	\begin{cl}\label{2510041}
		$F\nsubseteq A\cup B$ for any $F\in \left(\mf_{1}\bs \left\{B\right\}\right)\cup \left( \mf_{2}\bs \left\{A\right\}\right)$.
	\end{cl}
	\begin{proof}
		If $t=1$, then this claim is routine to check. In the following,  assume that $t\geq 2$.
		
		We first show $\left\{F\in\mf_{1}\bs\left\{B\right\}: F\nsubseteq A\cup B\right\}\neq \emptyset$.  Suppose for contradiction that $F\subseteq A\cup B$ for any $F\in\mf_{1}\bs\left\{B\right\}$. Then $F\subseteq A\cup B$ for any $F\in\mf_{1}\cup \mf_{3}\cup \left\{A\right\}$. Since $\left|\bigcup_{F\in\mf}F\right|\geq t+4$, there exits $F_{0}\in\mf_{2}\bs \left\{A\right\}$ such that $F_{0}\nsubseteq A\cup B$. Moreover, there are two members in $\mw(F_{0}, B)$ are  not in $\mf$, i.e., $\left| \mf\cap \mw(F_{0}, B)\right|\leq 2$. It follows from Proposition \ref{2509303} (i) and (iii) that 
		\begin{equation*}
			\begin{aligned}
				\left|\md_{\mf}(F_{0};t)\right|&\geq \left|\left\{B\right\}\cup\left(\mi_{\mf}(B;t)\bs\mw(F_{0}, B) \right)\right|=1+\left|\mi_{\mf}(B;t)\right|-\left|\mi_{\mf}(B;t)\cap \mw(F_{0}, B)\right|\\
				&\geq 1+\left|\mi_{\mf}(B;t)\right|-\left|\mf\cap \mw(F_{0}, B) \right|\geq s+1,
			\end{aligned}
		\end{equation*}
		a contradiction. Similarly, we have  $\left\{F\in\mf_{2}\bs\left\{A\right\}: F\nsubseteq A\cup B\right\}\neq \emptyset$.
		
		Now,  we prove this claim. Suppose for contradiction that $F_{0}\subseteq A\cup B$ for some $F_{0}\in \left(\mf_{1}\bs \left\{B\right\}\right)\cup \left( \mf_{2}\bs \left\{A\right\}\right)$. Then $A\cap B\nsubseteq F_{0}$. If $F_{0}\in \mf_{1}\bs \left\{B\right\}$, then pick $F_{1}\in \mf_{1}\bs\left\{B\right\}$ with $F_{1}\nsubseteq A\cup B$. By (\ref{2510042}), we have $A\cap B\subseteq B\cap F_{1}$. It follows from $A\cap B\nsubseteq F_{0}$ that
		$$F_{0}\cap F_{1}=F_{0}\cap \left(A\cup B\right)\cap F_{1}=F_{0}\cap B\cap F_{1}\subsetneq B\cap F_{1}.$$
		Hence $\left|F_{0}\cap F_{1}\right|\leq t-1$. By Proposition \ref{2510031}, we know	$\md_{\mf}(A;t)$ is $t$-intersecting, a contradiction to $F_{0}, F_{1}\in \md_{\mf}(A;t)$. If $F_{0}\in \mf_{2}\bs \left\{A\right\}$, then by the similar discussion, we also derive a contradiction.
	\end{proof}

From (\ref{2510042}) and Claim \ref{2510041}, we obtain
	\begin{equation}\label{2010043}
		\mf_{1}\bs \left\{B\right\}\subseteq  \mmu(A, B;b)\cup\mmu(A,B;d),\ \mf_{2}\bs \left\{A\right\}\subseteq   \mmu(A, B;a)\cup \mmu(A,B;c).
	\end{equation}
	By (\ref{2510033}), we know $\left|\mf_{1}\bs \left\{B \right\}\right|\geq 3$ and $\left|\mf_{2}\bs \left\{A\right\}\right|\geq 3$. Then  
	$\left| \mf_{1}\cap \mmu(A,B;i)\right|\geq 2$ and  $\left| \mf_{2}\cap \mmu(A,B;j)\right|\geq 2$ for some $i\in \left\{b,d\right\}$ and $j\in \left\{a,c\right\}$. 
	W.l.o.g, we  may assume that $i=b$ and $j=a$.
	\begin{cl}\label{2010044}
		$\mf_{1}\bs \left\{B\right\}\subseteq \mmu(A, B;b) $  and $\mf_{2}\bs  \left\{A\right\}\subseteq  \mmu(A, B;a)$. 
	\end{cl}
	\begin{proof}
		We only show 	$\mf_{1}\bs \left\{B\right\}\subseteq \mmu(A, B;b)$, the other proof is similar. Suppose for contradiction that $\mf_{1}\bs \left\{B\right\}\nsubseteq \mmu(A, B;b)$. It follows  from (\ref{2010043}) that there exists $m_{1}\in[n]\bs \left(A\cup B\right)$ such that $\left(A\cap B\right)\cup \left\{d,m_{1}\right\}\in \mf_{1}$. Since $\left| \mf_{1}\cap \mmu(A,B;b)\right|\geq 2$, we know $\left(A\cap B\right)\cup \left\{b,m_{2}\right\}\in \mf_{1}$ for some $m_{2}\in [n]\bs \left(A\cup B\right)$ with $m_{2}\neq m_{1}$. Then  $\left|\left(\left(A\cap B\right)\cup \left\{d,m_{1}\right\}\right)\cap \left(\left(A\cap B\right)\cup \left\{b,m_{2}\right\}\right) \right|=t-1$. By Proposition \ref{2510031}, we know	$\md_{\mf}(A;t)$ is $t$-intersecting. This contradicts to the fact that $\left(A\cap B\right)\cup \left\{d,m_{1}\right\}$ and $\left(A\cap B\right)\cup \left\{b,m_{2}\right\}$ are in $\md_{\mf}(A;t)$. 
	\end{proof}
	
	This claim together with $\mf_{3}\subseteq \mw(A, B)$ yields $2n-2t\geq \left|\mf_{1}\right|+\left|\mf_{2}\right|+\left|\mf_{3}\right|\geq 2s+3$, i.e., $n\geq t+s+2$.

	Applying Claim \ref{2010044} again, we know $\left(\mf_{1}\bs\left\{B\right\}\right)\cup \left(\mf_{2}\bs\left\{A\right\}\right)\subseteq \md_{\mf}(\left(A\cap B\right)\cup \left\{c,d\right\};t)$.
	It follows from (\ref{2510033}) that $\left|\md_{\mf}(\left(A\cap B\right)\cup \left\{c,d\right\};t)\right|\geq s+1$. Hence $\left(A\cap B\right)\cup \left\{c,d\right\}\notin \mf$. We further derive $\mf_{3}\subseteq \mw(A,B)\bs \left\{\left(A\cap B\right)\cup \left\{c,d\right\}\right\}$. This together with  $\left|\mf\right|\geq 2s+3$, $\left|\mf_{1}\right|\leq s$ and $\left|\mf_{2}\right|\leq s$ yields
	$$\mf_{3}= \mw(A,B)\bs \left\{\left(A\cap B\right)\cup \left\{c,d\right\}\right\},\ \left|\mf_{1}\bs \left\{B\right\}\right|=\left|\mf_{2}\bs \left\{A\right\}\right|=s-1.$$
	Pick $F\in \mf_{1}\bs \left\{B\right\}$. From Claim \ref{2010044}, there exists  $m\in [n]\bs\left(A\cup B\right)$ such that $F=\left(A\cap B\right)\cup \left\{b,m\right\}$.   Since $A$ and $\left(A\cap B\right)\cup\left\{a,d\right\}$ are in $\md_{\mf}(F;t)\bs \left( \mf_{2}\bs\left\{A\right\}\right)$, we know $\mf_{2}\bs \left\{ A\right\}\nsubseteq \md_{\mf}(F;t)$ due to $\left| \md_{\mf}(F;t)\right|\leq s$.  It follows from Claim \ref{2010044} that $\left(A\cap B\right)\cup \left\{a,m\right\}\in \mf_{2}\bs \left\{A\right\}$. Then there exists an $(s-1)$-subset $M$ of $[n]\bs\left(A\cup B\right)$ such that 
	$$\mf_{1}\bs \left\{B\right\}=\left\{ \left(A\cap B\right)\cup  \left\{b,m\right\}: m\in M\right\}, \ \mf_{2}\bs \left\{A\right\}=\left\{ \left(A\cap B\right)\cup \left\{a,m\right\}: m\in M\right\}.$$
	Therefore, we have $\mf=\mm_{1}(A\cup B\cup M, \left(A\cap B\right)\cup \left\{a,b\right\}, A\cap B;t)$.

	(ii) Since $\left|\bigcup_{F\in\mf}F\right|= t+3$, we have  $F\subseteq A\cup B$ for any $F\in\mf$, which implies
	\begin{equation}\label{2510034}
		\mf_{1}\bs \left\{B\right\}\subseteq  \mv(A, B;a)\cup\mv(A,B;c),\ \mf_{2}\bs \left\{A\right\}\subseteq   \mv(A, B;b)\cup \mv(A,B;d).
	\end{equation} 
	By (\ref{2510033}), we know $\left|\mf_{1}\bs \left\{ B\right\}\right|\geq 3$ and $\left|\mf_{2}\bs \left\{A\right\}\right|\geq 3$. Then   
	$\left| \mf_{1}\cap \mv(A,B;i)\right|\geq 2$ and  $\left| \mf_{2}\cap \mv(A,B;j)\right|\geq 2$ for some $i\in \left\{a,c\right\}$ and $j\in \left\{b,d\right\}$. 
	W.l.o.g, we may assume that $i=a$ and $j=b$.
	
	\begin{cl}\label{2510032}
		$\mf_{1}\bs \left\{B\right\}\subseteq \mv(A, B;a) $  and $\mf_{2}\bs  \left\{A\right\}\subseteq  \mv(A, B;b)$. 
	\end{cl}
	\begin{proof}
		We only show 	$\mf_{1}\bs \left\{B\right\}\subseteq \mv(A, B;a)$, the other proof is similar. Suppose for contradiction that $\mf_{1}\bs \left\{B\right\}\nsubseteq \mv(A, B;a)$. It follows from (\ref{2510034}) that $\left(A\cup B\right)\bs \left\{c,m_{1}\right\}\in\mf_{1}$ for some $m_{1}\in A\cap B$. Recall that $\left| \mf_{1}\cap \mv(A,B;a)\right|\geq 2$. There exists $m_{2}\in A\cap B$ such that $m_{2}\neq m_{1}$ and $\left(A\cup B\right)\bs \left\{a,m_{2}\right\}\in \mf_{1}$. 
		Note that $\left(A\cup B\right)\bs \left\{c,m_{1}\right\}\in\md_{\mf}(A;t)$ and 
		$\left(A\cup B\right)\bs \left\{a,m_{2}\right\}\in \md_{\mf}(A;t)\cap \md_{\mf}(\left(A\cup B\right)\bs \left\{c,m_{1}\right\};t)$.  
		By Proposition \ref{2510031}, we know	$\md_{\mf}(A;t)$ is $t$-intersecting, a contradiction.
	\end{proof}
	
	From this  claim and $\mf_{3}\subseteq \mw(A, B)$, we obtain $2t+4\geq \left|\mf_{1}\right|+\left|\mf_{2}\right|+\left|\mf_{3}\right|\geq 2s+3$, which implies $t\geq s$.
	
	Using Claim \ref{2510032} again, we have $\left(\mf_{1}\bs\left\{B\right\}\right)\cup \left(\mf_{2}\bs\left\{A\right\}\right)\subseteq \md_{\mf}(\left(A\cap B\right)\cup \left\{a,b\right\};t)$.
	This together  with (\ref{2510033}) yields $\left|\md_{\mf}(\left(A\cap B\right)\cup \left\{a,b\right\};t)\right|\geq s+1$, which implies $\left(A\cap B\right)\cup \left\{a,b\right\}\notin \mf$.   Then $\mf_{3}\subseteq \mw(A,B)\bs \left\{\left(A\cap B\right)\cup \left\{a,b\right\}\right\}$. It follows from $\left|\mf\right|\geq 2s+3$, $\left|\mf_{1}\right|\leq s$ and $\left|\mf_{2}\right|\leq s$ that 
	$$\mf_{3}= \mw(A,B)\bs \left\{\left(A\cap B\right)\cup \left\{a,b\right\}\right\},\ \left|\mf_{1}\bs\left\{B\right\}\right|=\left|\mf_{2}\bs \left\{A\right\}\right|=s-1.$$
	Let $F\in\mf_{1}\bs \left\{B\right\}$. By  Claim \ref{2510032}, we get $F=\left(A\cup B\right)\bs \left\{ a,m\right\}$ for some $m\in A\cap B$. Note that $A$ and $\left(A\cap B\right)\cup \left\{ a,d\right\}$ are in $\md_{\mf}(F;t)\bs \left(\mf_{2}\bs\left\{A\right\}\right)$. Then $\mf_{2}\bs \left\{A\right\} \nsubseteq \md_{\mf}(F;t)$. This combining with  Claim \ref{2510032} that $\left( A\cup B\right)\bs \left\{b, m\right\}\in \mf_{2}\bs \left\{A\right\}$. Hence there exists an $(s-1)$-subset $M$ of $A\cap B$ such that 
	$$\mf_{1}\bs \left\{B\right\}=\left\{ \left(A\cup B\right)\bs \left\{a,m\right\}: m\in M\right\}, \ \mf_{2}\bs \left\{A\right\}=\left\{ \left(A\cup B\right)\bs \left\{b,m\right\}: m\in M\right\}. $$ 
	Therefore, we know $\mf=\mm_{2}(A\cup B, \left\{a,b\right\}, \left(A\cap B\right)\bs M;t)$.
\end{proof}

\begin{proof}[\bf Proof of Theorem \ref{2509213}]
For each   $s$-almost $t$-intersecting family  $\mf\subseteq \binom{[n]}{t+1}$,  if it is  not $t$-intersecting, then $\left|\bigcup_{F\in\mf}F\right|\geq t+3$. The desired result  follows from Lemma \ref{2510311}.
\end{proof}

\begin{re}\label{2601021}
 The condition $s \geq 7$ in Theorem \ref{2509213}  cannot be omitted.   In fact, for $s \in \{1,3,6\}$, there exist  families  distinct from the extremal families in Theorem \ref{2509213} whose  sizes are at least $2s + 3$.

Let   $X$ and $Y$ be two subsets of $[n]$ with $\left|X\right|\geq t+3$,  $\left|Y\right|\leq t-1$ and $Y\subseteq X$.
Set 	$\mf=\left\{ F\in \binom{X}{t+1}: Y\subseteq F\right\}$, which  is not $t$-intersecting. Moreover, 
 if $\left( \left|X\right|, \left|Y\right|\right)=(t+3,t-1)$, then $\mf$ is $1$-almost $t$-intersecting; 	if $\left( \left|X\right|,\left|Y\right|\right)=(t+4,t-1)$, or $t\geq 2$ and  $\left( \left|X\right|,\left|Y\right|\right)=(t+3,t-2)$, then $\mf$ is $3$-almost $t$-intersecting;
	if $\left( \left|X\right|,\left|Y\right|\right)=(t+5,t-1)$, or $t\geq 3$ and  $\left( \left|X\right|,\left|Y\right|\right)=(t+3,t-3)$, then $\mf$ is $6$-almost $t$-intersecting. 
\end{re}

\section{Inequalities concerning binomial coefficients}\label{2510202}

This section establishes some inequalities for binomial coefficients, beginning with a lemma proved via routine computations.

\begin{lem}\label{2510071}
  Let $n$, $k$, $t$ and $s$ be positive integers with $k\geq t+2$. The following hold.
  \begin{enumerate}[\normalfont(i)]
  	\item If $n\geq L_{1}(k,t,s)$, then
  	$$s\binom{2k-2t+2}{k-t+1}<\min\left\{ \binom{n-k-1}{k-t}-t, \frac{3(k-t+1)s^{\frac{1}{k-t}}}{k-t}\binom{n-t-2}{k-t-1}\right\}.$$
  	\item If $n\geq L_{2}(k,t,s)$,   then 
  	\begin{equation*}
  		\frac{11}{54}(k-t+1)\binom{n-t-1}{k-t-1}\geq \frac{33}{54}\binom{n-t-1}{k-t-1}>\left\{
  		\begin{aligned}
  			&11s+9t+18 &\textnormal{if}\ \  k=t+2,\\
  			&s\binom{2k-2t+2}{k-t+1} &\textnormal{if}\ \  k\geq t+3.
  		\end{aligned}
  		\right.
  	\end{equation*}
  \end{enumerate}
\end{lem}
\begin{proof}
 (i)	Notice that $(t+1)(k-t+1)\geq 2k$. For each $i\in \left\{0,1,\ldots, k-t-1\right\}$, we have
 $$n-k-i-1\geq 3(k-t+1)s^{\frac{1}{k-t}}+k-i-1\geq \frac{3s^{\frac{1}{k-t}}}{2}(2k-2t-i+1)+t. $$
 It follows that 
 \begin{equation*}
 	\begin{aligned} 
 		\binom{n-k-1}{k-t}&=\prod_{i=0}^{k-t-1}\frac{n-k-i-1}{k-t-i}\geq \left(\frac{3s^{\frac{1}{k-t}}(k-t+2)}{2}+t\right)\prod_{i=0}^{k-t-2}\frac{3s^{\frac{1}{k-t}}(2k-2t-i+1)}{2(k-t-i)}\\
 		&= s\left(\frac{3}{2}\right)^{k-t}\binom{2k-2t+1}{k-t}+t\prod_{i=0}^{k-t-2}\frac{3s^{\frac{1}{k-t}}(2k-2t-i+1)}{2(k-t-i)}\\
 		&\geq \frac{9s}{4}\binom{2k-2t+1}{k-t}+t>s\binom{2k-2t+2}{k-t+1}+t,
 	\end{aligned}
 \end{equation*}
 \begin{equation*}
 	\begin{aligned}
 		\frac{3(k-t+1)s^{\frac{1}{k-t}}}{k-t}\binom{n-t-2}{k-t-1}&=\frac{3(k-t+1)s^{\frac{1}{k-t}}}{k-t}\prod_{i=1}^{k-t-1}\frac{n-t-i-1}{k-t-i}\\
 		&\geq \frac{3(2k-2t+1)s^{\frac{1}{k-t}}}{2(k-t)}\prod_{i=1}^{k-t-1}\frac{3s^{\frac{1}{k-t}}(2k-2t-i+1)}{2(k-t-i)}\\
 		&= s\left(\frac{3}{2}\right)^{k-t}\binom{2k-2t+1}{k-t}>s\binom{2k-2t+2}{k-t+1}.
 	\end{aligned}
 \end{equation*}
 Then (i) holds.
 
 (ii) If $k=t+2$, then we have 
 \begin{equation*}
 	\begin{aligned}
 		3(n-t-1)\geq 54s+105t+54t^{2}+51\geq 54s+105t+105> \frac{54}{11}\left(11s+9t+18\right). 
 	\end{aligned}
 \end{equation*}
 If $k\geq t+3$, then for each $i\in \left\{0,1,\ldots, k-t-2\right\}$,
 $$n-t-i-1\geq 6(k-t+1)s^{\frac{1}{k-t-1}}-i\geq 3(2k-2t-i)s^{\frac{1}{k-t-1}}. $$
 We further derive 
 \begin{equation*}
 	\begin{aligned}
 		3\binom{n-t-1}{k-t-1}&=3\prod_{i=0}^{k-t-2}\frac{n-t-i-1}{k-t-i-1}\geq 3\prod_{i=0}^{k-t-2}\frac{3(2k-2t-i)s^{\frac{1}{k-t-1}}}{k-t-i-1}\\
 		&= s\cdot \frac{3^{k-t}}{2}\left(\frac{1}{2}-\frac{1}{2(2k-2t+1)} \right)\binom{2k-2t+2}{k-t+1}\\
 		&\geq \frac{81}{14}s\binom{2k-2t+2}{k-t+1}>\frac{54}{11}s\binom{2k-2t+2}{k-t+1}.
 	\end{aligned}
 \end{equation*}
Note that $k-t+1\geq 3$. Hence (ii) follows.  
\end{proof}

Recall that $f_{1}(n,k,t,x,y)$, $f_{2}(n,k,t,s,x)$, $f_{3}(n,k,t,s)$  are defined in (\ref{2510082}), (\ref{2509295}), (\ref{2509296}) respectively, and $g_{1}(n,k,t)$, $g_{2}(n,k,t)$, $g_{3}(n,k,t,s)$ are defined in (\ref{2510162}), (\ref{2510163}), (\ref{2509244}) respectively. The following lemmas show some inequalities concerning these functions.

\begin{lem}\label{2510072}
	 Let $n$, $k$, $t$ and $s$ be positive integers with  $n\geq L_{1}(k,t,s)$. The following hold.
	 \begin{enumerate}[\normalfont(i)]
	 	\item If $k=t+1$, then $n-t>\max\left\{t+2, 2s+4\right\}$.
	 	\item If $k\geq t+2$, then $\binom{n-t}{k-t}> \max\left\{ g_{1}(n,k,t),g_{2}(n,k,t)\right\}+s\binom{2k-2t+2}{k-t+1}$.
	 \end{enumerate}
\end{lem}
\begin{proof}
	(i) From  $n-t\geq t+6s+2>\max\left\{ t+2, 2s+4\right\}$, we obtain (i).

	(ii) 
    Using Pascal’s formula, we get 
$ \binom{n-t}{k-t}-g_{1}(n,k,t)= \binom{n-t-1}{k-t}-(t+1)\binom{n-t-2}{k-t-1}$.
	It follows Lemma \ref{2510071} (i) that 
	\begin{equation*}
		\begin{aligned}
					&\  \binom{n-t}{k-t}-g_{1}(n,k,t)-s\binom{2k-2t+2}{k-t+1}\\
			\geq&\ \frac{(k-t+1)((t+1)+3s^{\frac{1}{k-t}})-(t+1)}{k-t}\binom{n-t-2}{k-t-1}-(t+1)\binom{n-t-2}{k-t-1}-s\binom{2k-2t+2}{k-t+1}\\
			=& \frac{3(k-t+1)s^{\frac{1}{k-t}}}{k-t}\binom{n-t-2}{k-t-1}-s\binom{2k-2t+2}{k-t+1}>0.
		\end{aligned}
	\end{equation*}
 Applying  Lemma \ref{2510071} (i) again, we get $\binom{n-t}{k-t}-g_{2}(n,k,t)-s\binom{ 2k-2t+2}{k-t+1}>0$.
\end{proof}

\begin{pr1}\label{2510083}\textnormal{(\cite[Proposition 1.6]{2504063})}
	Let $n$, $k$ and $i$ be positive integers with $n>ik$. Then $\binom{n-i}{k}\geq \left( 1-\frac{ik}{n}\right) \binom{n}{k}$.
\end{pr1}

In the following, we always assume that $n,k,t$ and $s$ are positive integers with $k\geq t+2$ and $n\geq L_{2}(k,t,s)$. Then 
\begin{equation}\label{2510084}
	n-t-2\geq \max\left\{2(t+1)^{2}(k-t)(k-t+1), 6(k-t+1)s^{\frac{1}{k-t-1}} \right\}.
\end{equation}

\begin{lem}\label{2510081}
	$f_{1}(n,k,t,t+2,t+2)\geq f_{1}(n,k,t,x,x)$ for any $x\in\left\{t+2, t+3, \ldots, k\right\}$.
\end{lem}
\begin{proof}
	If $k=t+2$, then there is nothing to prove. Next, we assume that $k\geq t+3$. Let $x\in\left\{t+2, t+3,\ldots, k-1\right\}$. Note that $\left(1+\frac{1}{k-2t+x}\right)^{x-t}$ is no more than  Euler's number $e$. It follows from (\ref{2510084}) that
	\begin{equation*}
		\begin{aligned}
			\frac{f_{1}(n,k,t,x,x)}{f_{1}(n,k,t,x+1,x+1)}&=\frac{x-t+1}{x+1}\cdot  \frac{1}{\left(1+\frac{1}{k-2t+x}\right)^{x-t}}\cdot \frac{1}{k-2t+x+1}\cdot \frac{n-x}{k-x}\\
			&\geq \frac{3}{t+3}\cdot \frac{1}{e}\cdot \frac{1}{2(k-t)}\cdot \frac{n-t-2}{k-t-2}\geq \frac{3}{t+3}\cdot \frac{1}{e}\cdot (t+1)^{2}\\
			&=\frac{3(t+1)}{e}\cdot\frac{t+1}{t+3}\geq 1,
		\end{aligned}
	\end{equation*}
	which implies the desired result.
\end{proof}

\begin{lem}\label{2510114}
	$f_{2}(n,k,t,s,k-t)\geq f_{2}(n,k,t,s,x)$ for any $x\in\left\{0, 1, \ldots, k-t\right\}$.
\end{lem}
\begin{proof}
	For each $x\in\left\{0, 1, \ldots, k-t-1\right\}$,  we have 
	\begin{equation}\label{2510085}
		\begin{aligned}
			&f_{2}(n,k,t,s,x+1)-f_{2}(n,k,t,s,x)\\
			=&\binom{n-t-x-1}{k-t-1}-\left(2k-2t-2x+1\right)\binom{n-t-2}{k-t-2}\\
			\geq& \binom{n-k}{k-t-1}-\left(2k-2t+1\right)\binom{n-t-2}{k-t-2}.
		\end{aligned}
	\end{equation}
	By Proposition \ref{2510083}, we know $\binom{n-k}{k-t-1}\geq \left(1-\frac{(k-t-1)^{2}}{n-t-1}\right)\binom{n-t-1}{k-t-1}$. It follows from (\ref{2510084}) and (\ref{2510085}) that 	for each $x\in\left\{0, 1, \ldots, k-t-1\right\}$,
	\begin{equation*}
		\begin{aligned}
			&f_{2}(n,k,t,s,x+1)-f_{2}(n,k,t,s,x)\\
			\geq & \left(1-\frac{(k-t-1)^{2}}{n-t-1}\right)\binom{n-t-1}{k-t-1}-\left(2k-2t+1\right)\binom{n-t-2}{k-t-2}\\
			=& \left( \frac{n-t-1}{k-t-1}-3(k-t)\right)\binom{n-t-2}{k-t-2}\\
			\geq&\left( 2(t+1)^{2}(k-t+1)-3(k-t)\right)\binom{n-t-2}{k-t-2}>0.
		\end{aligned}
	\end{equation*}
This implies the required result.
\end{proof}

\begin{lem}\label{2510086}
	The following hold.
	\begin{enumerate}[\normalfont(i)]
		\item  $ g_{3}(n,k,t,s)>f_{1}(n,k,t,t+2,t+2) +\frac{11}{54}(k-t+1)\binom{n-t-1}{k-t-1}$.
		\item $ g_{3}(n,k,t,s)>f_{1}(n,k,t,t+1,t+2)+\left(\frac{11}{54}(k-t+1)+1\right) \binom{n-t-1}{k-t-1}$.
		\item $ g_{3}(n,k,t,s)>f_{3}(n,k,t,s)$. 
	\end{enumerate}
\end{lem}
\begin{proof}
	By applying Pascal’s formula repeatedly, we get
	\begin{equation*}
		\binom{n-t}{k-t}-\binom{n-k-1}{k-t}=\sum_{i=0}^{k-t}\binom{n-t-i-1}{k-t-1}.
	\end{equation*}
	From Proposition \ref{2510083}, for each $i\in \left\{1,2,\ldots, k-t\right\}$,  we obtain $$\binom{n-t-i-1}{k-t-1}\geq \left( 1-\frac{i(k-t-1)}{n-t-1}\right)\binom{n-t-1}{k-t-1}.$$ Therefore, we have
	\begin{equation*}
		\begin{aligned}
			\frac{ 	\binom{n-t}{k-t}-\binom{n-k-1}{k-t}}{(k-t+1)\binom{n-t-1}{k-t-1}}\geq  \frac{1}{k-t+1}\sum_{i=0}^{k-t}\left(1-\frac{i(k-t-1)}{n-t-1}\right)= 1-\frac{(k-t-1)(k-t)}{2(n-t-1)}.
		\end{aligned}
	\end{equation*}
 	It follows from (\ref{2510084})   that 
	\begin{equation}\label{2510101}
		 g_{3}(n,k,t,s)>\binom{n-t}{k-t}-\binom{n-k-1}{k-t}\geq \frac{15}{16} (k-t+1)\binom{n-t-1}{k-t-1}.
	\end{equation}

(i) 	By applying (\ref{2510084}), we know 
	\begin{equation*}
		\begin{aligned}
			\frac{f_{1}(n,k,t,t+2,t+2)}{(k-t+1)\binom{n-t-1}{k-t-1}}&=\frac{(k-t-1)(k-t+2)^{2}\binom{t+2}{2}}{(k-t+1)(n-t-1)}\leq \frac{ (k-t-1)(k-t+2)^{2}\binom{t+2}{2}}{2(k-t)(k-t+1)^{2}(t+1)^{2}}\\
			&=\frac{1}{4}\cdot \frac{k-t-1}{k-t} \left( 1+\frac{1}{t+1}\right)\left(1+\frac{1}{k-t+1}\right)^{2}\leq \frac{2}{3}.
		\end{aligned}
	\end{equation*}
This together with  (\ref{2510101}) yields (i). 	

  (ii) From (\ref{2510084}), we get
     \begin{equation*}
     	\begin{aligned}
     	\frac{11}{54}+	\frac{\binom{n-t-1}{k-t-1}+f_{1}(n,k,t,t+1,t+2)}{(k-t+1)\binom{n-t-1}{k-t-1}}&=\frac{11}{54}+\frac{1}{k-t+1}+\frac{(t+1)(k-t-1)(k-t+1)}{n-t-1}\\
     		&\leq\frac{11}{54}+ \frac{1}{k-t+1}+\frac{k-t-1}{2(t+1)(k-t)}\leq \frac{85}{108}.
     	\end{aligned}    	
     \end{equation*}
It follows from (\ref{2510101}) that (ii) holds.

 (iii) 	By (\ref{2510084}) and Lemma \ref{2510071} (ii), we have
  \begin{equation}\label{2510112}
 	\frac{n-t-1}{2(k-t-1)}\geq (t+1)(k-t)+1,\ \  \frac{1}{2}\binom{n-t-1}{k-t-1}>2s.
  \end{equation}

	From Proposition \ref{2510083}, we obtain $\binom{n-k-1}{k-t-1}\geq \left(1-\frac{(k-t)(k-t-1)}{n-t-1}\right)\binom{n-t-1}{k-t-1}$. It follows  that
	\begin{equation*}
		\begin{aligned}
			& g_{3}(n,k,t,s)-f_{3}(n,k,t,s)\\
			=&\binom{n-k-1}{k-t-1}-\left(t(k-t)+1\right)\binom{n-t-2}{k-t-2}-2s+\min\left\{t,s\right\}\\
			\geq & \left(1-\frac{(k-t)(k-t-1)}{n-t-1}\right)\binom{n-t-1}{k-t-1}-\left(t(k-t)+1\right)\binom{n-t-2}{k-t-2}-2s\\
			=& \left( \frac{n-t-1}{2(k-t-1)}-(t+1)(k-t)-1\right)\binom{n-t-2}{k-t-2}+\frac{1}{2}\binom{n-t-1}{k-t-1}-2s.
		\end{aligned}
	\end{equation*}
	This combining with (\ref{2510112}) yields (iii). 
\end{proof}

\medskip
\noindent{\bf Acknowledgment.}	

The authors would like to thank Professor Alex Scott for his  valuable comments.	K. Wang \\is supported by the National Natural Science Foundation of China (12131011, 12571347) and Beijing Natural Science Foundation (1252010, 1262010). T. Yao is supported by   Natural Science Foundation of Henan (252300420899).

\end{document}